\documentclass{article}

\input{tcilatex}
\begin{document}

\section{1. Introduction}

\qquad Let $p(x)$ be a polynomial of degree $n\geq 2$ with $n$ distinct real
roots $r_{1}<r_{2}<\cdots <r_{n}$. Such a polynomial is called hyperbolic.
Let $x_{1}<x_{2}<\cdots <x_{n-1}$ be the critical points of $p$, and define
the ratios $\sigma _{k}=\dfrac{x_{k}-r_{k}}{r_{k+1}-r_{k}},k=1,2,...,n-1$. $%
(\sigma _{1},...,\sigma _{n-1})$ is called the \textit{ratio vector} of $p$,
and $\sigma _{k}$ is called the $k$th ratio. Ratio vectors were first
discussed in [5] and in [1], where the inequality $\dfrac{1}{n-k+1}<\sigma
_{k}<\dfrac{k}{k+1},k=1,2,...,n-1$ was derived. In a similar fashion, one
can define ratios for polynomial like functions of the form $%
p(x)=(x-r_{1})^{m_{1}}\cdots (x-r_{N})^{m_{N}}$, where $m_{1},...,m_{N}$ are
given positive real numbers and $r_{1}<r_{2}<\cdots <r_{N}$(see [4]).

In this paper we want to discuss the extension of the notion of ratios to
polynomials with \textit{complex} roots. Thus we let $p(z)$ be a polynomial
of degree $n\geq 2$ with $n$ distinct complex roots $w_{1},...,w_{n}$ and
critical points $z_{1},...,z_{n-1}$. Numerous papers have investigated the
relation between the roots and critical points of a polynomial. The focus of
this paper is to investigate that relation in the form of the complex ratios 
$\sigma _{k}=\dfrac{z_{k}-w_{k}}{w_{k+1}-w_{k}},k=1,2,...,n-1$. The main
problem is in defining the ratios when there is no natural ordering of roots
and critical points as with all real roots. We have to order the $\left\{
w_{k}\right\} $ somehow and then determine which $\left\{ z_{k}\right\} $
are associated with $w_{k}$ and $w_{k+1}$. We use the real parts of the $%
\left\{ w_{k}\right\} $ and the $\left\{ z_{k}\right\} $ to do this. For the
rest of the paper we concentrate solely on the case $n=3$, which is already
fairly nontrivial. We do not define the ratios in the case when two roots or
critical points have equal real parts(unless the critical points are
identical). One could certainly extend the definition to those cases, but
the ratios will not be continuous function of the roots. Our definition does
extend the definition of the ratios when $p$ is hyperbolic and the ratios
are continuous functions of the roots when the roots are all real. For cubic
hyperbolic polynomials, the inequality $\dfrac{1}{n-k+1}<\sigma _{k}<\dfrac{k%
}{k+1}$ implies that $\dfrac{1}{3}<\sigma _{1}<\dfrac{1}{2}\,$and $\dfrac{1}{%
2}<\sigma _{2}<\dfrac{2}{3}$. For complex ratios, we derive separate and
sharp upper and lower bounds on the real and imaginary parts, and modulus,
of each ratio(see Theorems 1 and 2). For cubic hyperbolic polynomials, it is
immediate that $\sigma _{1}<\sigma _{2}$. In the complex case we prove that $%
\func{Re}\sigma _{1}\leq \func{Re}\sigma _{2}$(Theorem 3). Indeed, one can
have $\sigma _{1}=\sigma _{2}$(see Theorem 4). Finally, we show that the
ratios are real if and only if the roots of $p$ are collinear(Theorem 5).

\section{2. Main Results}

Let 
\[
p(w)=(w-w_{1})(w-w_{2})(w-w_{3}), 
\]%
where we assume that $\func{Re}w_{1}<\func{Re}w_{2}<\func{Re}w_{3}$. Now the
critical points of $p$ are $\dfrac{1}{3}\left( w_{1}+w_{2}+w_{3}\pm \sqrt{%
w_{1}^{2}+w_{2}^{2}+w_{3}^{2}-w_{1}w_{3}-w_{1}w_{2}-w_{2}w_{3}}\right) $,
where $\sqrt{z}$ is the principal branch of the square root function,
analytic everywhere except on the \textit{nonpositive} real axis, which we
denote by $\Gamma $. Note that $\func{Re}\sqrt{z}\geq 0$ and $\sqrt{z^{2}}=z$
if $\func{Re}z\geq 0$. We also assume that if the critical points are not
identical, then they cannot have equal real parts. In other words, we assume
that

$\func{Re}\sqrt{%
w_{1}^{2}+w_{2}^{2}+w_{3}^{2}-w_{1}w_{3}-w_{1}w_{2}-w_{2}w_{3}}\neq 0$
unless $w_{1}^{2}+w_{2}^{2}+w_{3}^{2}-w_{1}w_{3}-w_{1}w_{2}-w_{2}w_{3}=0$.
Denote the critical points by $z_{1}$ and $z_{2}$, where $z_{1}=z_{2}$, or $%
\func{Re}z_{1}<\func{Re}z_{2}$ if $z_{1}\neq z_{2}$. We define the ratios

\begin{equation}
\sigma _{1}=\dfrac{z_{1}-w_{1}}{w_{2}-w_{1}},\sigma _{2}=\dfrac{z_{2}-w_{2}}{%
w_{3}-w_{2}}  \tag{(1)}
\end{equation}

$(\sigma _{1},\sigma _{2})$ is called the \textit{ratio vector} of $p$. One
can give a geometric interpretation for the ratios as follows. First, if $%
p\in \pi _{3}$ has \textit{noncollinear} zeros, let $T$ be the triangle
whose vertices are $w_{1},w_{2},w_{3}$. Let $E$ be the midpoint ellipse,
that is, the ellipse tangent to $T$ at the midpoints of its sides. Then it
is well known that the zeros of $p^{\prime }$ are the foci, $z_{1}$ and $%
z_{2}$, of $E$. Let $\theta _{1}$ denote the angle between $\overrightarrow{%
w_{1}z_{1}}$ and $\overrightarrow{w_{1}w_{2}}$, and let $\theta _{2}$ denote
the angle between $\overrightarrow{w_{2}z_{2}}$ and $\overrightarrow{%
w_{2}w_{3}}$. Then $\theta _{1}=\arg \sigma _{1}$ and $\theta _{2}=\arg
\sigma _{2}$, and thus each ratio represents an angle between a line segment
of the circumscribed triangle of $E$ and a line segment connecting a vertex
to one of the foci of $E$.

Using our definition, neither $\sigma _{1}$ nor $\sigma _{2}$ will be a
continuous function of $w_{1},w_{2},$ and $w_{3}$, though they are
continuous on an open subset of $C^{3}-\left\{
w_{1},w_{2},w_{3}:w_{i}=w_{j}\ \text{for some }i\neq j\right\} $ and in
particular at any point $(w_{1},w_{2},w_{3})$ where all of the $\left\{
w_{k}\right\} $ are real. Clearly, if we translate the roots of $p$, the
ratios $\sigma _{1}$ and $\sigma _{2}$ do not change. Thus we may assume
that 
\begin{equation}
w_{1}+w_{3}=0  \tag{(2)}
\end{equation}

which implies that $\func{Re}w_{1}<0<\func{Re}w_{3}$. Note that $\func{Re}%
\sqrt{w_{3}}>0$. The critical points of $p$ are then $\dfrac{1}{3}\left(
w_{2}\pm \sqrt{3w_{3}^{2}+w_{2}^{2}}\right) $. The assumption that if the
critical points are not identical, then they cannot have equal real parts
now takes the form 
\[
3w_{3}^{2}+w_{2}^{2}\neq 0\Rightarrow \func{Re}\sqrt{3w_{3}^{2}+w_{2}^{2}}%
\neq 0 
\]

If $3w_{3}^{2}+w_{2}^{2}\neq 0$, then by our choice of the branch of $\sqrt{z%
},\func{Re}\sqrt{3w_{3}^{2}+w_{2}^{2}}>0$, which implies that $\func{Re}%
\left( w_{2}-\sqrt{3w_{3}^{2}+w_{2}^{2}}\right) <\func{Re}\left( w_{2}+\sqrt{%
3w_{3}^{2}+w_{2}^{2}}\right) $. Thus we have

\[
z_{1}=\dfrac{1}{3}\left( w_{2}-\sqrt{3w_{3}^{2}+w_{2}^{2}}\right) \text{, }%
z_{2}=\dfrac{1}{3}\left( w_{2}+\sqrt{3w_{3}^{2}+w_{2}^{2}}\right) 
\]

Also, $\func{Re}\sqrt{3w_{3}^{2}+w_{2}^{2}}>0\iff 3w_{3}^{2}+w_{2}^{2}\notin
\Gamma $. That leads to the following.

\textbf{Definition: }We say that $(w_{2},w_{3})$ is an admissible pair if $%
w_{2}$ and $w_{3}$ satisfy $3w_{3}^{2}+w_{2}^{2}\notin \Gamma $, $%
w_{2}+w_{3}\neq 0$, $\func{Re}w_{2}<\func{Re}w_{3},$ and $0<\func{Re}w_{3}$.
A region in $C^{2}$ consisting of only admissible pairs is also called
admissible.

Note that the ratios are not defined, say, when $w_{1}=-1,w_{2}=ti,w_{3}=1$, 
$\left\vert t\right\vert >\sqrt{3}$, since in that case $\func{Re}\sqrt{%
3w_{3}^{2}+w_{2}^{2}}=0$, which implies that $\func{Re}z_{1}=\func{Re}z_{2}$%
, but $z_{1}\neq z_{2}$. Let 
\begin{equation}
w=\dfrac{w_{2}}{w_{3}}.  \tag{(3)}
\end{equation}%
We shall express $\sigma _{1}$ and $\sigma _{2}$ as analytic functions of $w$%
. We then derive bounds on the real part, imaginary part, and modulus of the
ratios and also some relations between the ratios. By (1) and (2), $\sigma
_{1}=\dfrac{z_{1}+w_{3}}{w_{2}+w_{3}}=\dfrac{1}{3}\dfrac{3z_{1}+3w_{3}}{%
w_{2}+w_{3}}=\dfrac{1}{3}\dfrac{w_{2}+3w_{3}-\sqrt{3w_{3}^{2}+w_{2}^{2}}}{%
w_{2}+w_{3}}=\dfrac{1}{3}\dfrac{w_{2}+3w_{3}-\sqrt{w_{3}^{2}\left( 3+\dfrac{%
w_{2}^{2}}{w_{3}^{2}}\right) }}{w_{2}+w_{3}}$. In general, $\sqrt{%
w_{3}^{2}\left( 3+\dfrac{w_{2}^{2}}{w_{3}^{2}}\right) }=\pm \sqrt{w_{3}^{2}}%
\sqrt{3+\dfrac{w_{2}^{2}}{w_{3}^{2}}}=w_{3}\sqrt{3+\dfrac{w_{2}^{2}}{%
w_{3}^{2}}}$ or $-w_{3}\sqrt{3+\dfrac{w_{2}^{2}}{w_{3}^{2}}}$and thus $%
\sigma _{1}=f_{1}(w_{2},w_{3})$ or $\sigma _{1}=f_{2}(w_{2},w_{3})$, where $%
f_{1}(w_{2},w_{3})=\dfrac{1}{3}\dfrac{\dfrac{w_{2}}{w_{3}}+3-\sqrt{3+\dfrac{%
w_{2}^{2}}{w_{3}^{2}}}}{\dfrac{w_{2}}{w_{3}}+1}$ and $f_{2}(w_{2},w_{3})=%
\dfrac{1}{3}\dfrac{\dfrac{w_{2}}{w_{3}}+3+\sqrt{3+\dfrac{w_{2}^{2}}{w_{3}^{2}%
}}}{\dfrac{w_{2}}{w_{3}}+1}$. Now $\dfrac{1}{3}\dfrac{w_{2}+3w_{3}-\sqrt{%
3w_{3}^{2}+w_{2}^{2}}}{w_{2}+w_{3}}$ must be an analytic function of $w_{2}$
and $w_{3}$ in any admissible region. $f_{1}$ and $f_{2}$ are also analytic
functions of $w_{2}$ and $w_{3}$ in any admissible region with the
additional assumption that 
\begin{equation}
3+\dfrac{w_{2}^{2}}{w_{3}^{2}}\notin \Gamma .  \tag{(4)}
\end{equation}%
Since $3+\dfrac{w_{2}^{2}}{w_{3}^{2}}\neq 0,$ it follows that $\dfrac{1}{3}%
\dfrac{w_{2}+3w_{3}-\sqrt{3w_{3}^{2}+w_{2}^{2}}}{w_{2}+w_{3}}$ must equal $%
f_{1}(w_{2},w_{3})$ or $f_{2}(w_{2},w_{3})$. Now $f_{1}(0,1)=\allowbreak 1-%
\dfrac{1}{3}\sqrt{3}$, while $f_{2}(0,1)=\allowbreak 1+\dfrac{1}{3}\sqrt{3}$%
. But $w_{2}=0$ and $w_{3}=1$ yields the polynomial $p(z)=z(z^{2}-1)$, and
it is easy to check that $\sigma _{1}=1-\dfrac{1}{3}\sqrt{3}$. It then
follows that $\dfrac{1}{3}\dfrac{w_{2}+3w_{3}-\sqrt{3w_{3}^{2}+w_{2}^{2}}}{%
w_{2}+w_{3}}=\dfrac{1}{3}\dfrac{\dfrac{w_{2}}{w_{3}}+3+\sqrt{3+\dfrac{%
w_{2}^{2}}{w_{3}^{2}}}}{\dfrac{w_{2}}{w_{3}}+1}$. Using (3) we have $\sigma
_{1}=\dfrac{1}{3}\dfrac{w+3-\sqrt{3+w^{2}}}{w+1}$. Now let 
\[
E=\left\{ w:\func{Re}w=0,\left\vert \func{Im}w\right\vert \geq \sqrt{3}%
\right\} 
\]%
and let 
\[
D_{1}=C^{2}-E-\left\{ w:w=-1\right\} ,D_{2}=C^{2}-E-\left\{ w:w=1\right\} . 
\]

Note that $w\in C^{2}-E\iff (w_{2},w_{3})$ satisfies (4). Then 
\[
\sigma _{1}=\dfrac{1}{3}\dfrac{w+3-\sqrt{3+w^{2}}}{w+1},w\in D_{1}. 
\]%
In a similar fashion one can show that 
\[
\sigma _{2}=\dfrac{1}{3}\dfrac{-2w+\sqrt{3+w^{2}}}{1-w},w\in D_{2}. 
\]%
This expression for $\sigma _{2}$ also follows from the equation 
\begin{equation}
\left( 1-\sigma _{1}\right) \sigma _{2}=\dfrac{1}{3}  \tag{(5)}
\end{equation}

(5) is easy to prove and the proof is exactly the same as for the case when $%
p$ has three distinct real roots (see [1] or [3]). It is now convenient to
define the following analytic extensions of $\sigma _{1}$ to $w=-1$ and of $%
\sigma _{2}$ to $w=1$, respectively.

\[
f(w)=\left\{ 
\begin{array}{ll}
\dfrac{1}{3}\dfrac{w+3-\sqrt{3+w^{2}}}{w+1}=\sigma _{1} & \text{if }w\in
D_{1} \\ 
\dfrac{1}{2} & \text{if }w=-1%
\end{array}%
\right. 
\]%
and 
\[
g(w)=\left\{ 
\begin{array}{ll}
\dfrac{1}{3}\dfrac{-2w+\sqrt{3+w^{2}}}{1-w}=\sigma _{2} & \text{if }w\in
D_{2} \\ 
\dfrac{1}{2} & \text{if }w=1%
\end{array}%
\right. 
\]

Since $\lim\limits_{w\rightarrow -1}\dfrac{1}{3}\dfrac{w+3-\sqrt{3+w^{2}}}{%
w+1}=\dfrac{1}{2}$ and $\lim\limits_{w\rightarrow 1}\dfrac{1}{3}\dfrac{-2w+%
\sqrt{3+w^{2}}}{1-w}=\dfrac{1}{2}$, $f$ and $g$ are each analytic in the
region%
\[
D=C^{2}-E. 
\]

We can now replace (5) by

\begin{equation}
(1-f(w))g(w)=\dfrac{1}{3},w\in D  \tag{(6)}
\end{equation}

Note that $f$ does \textit{not} extend to be continuous on $\partial \left(
D\right) $ because of the discontinuity of $\sqrt{3+w^{2}}$ when $3+w^{2}\in
\Gamma $. Also, for $w\in \partial \left( D\right) ,f(w)$ does not yield $%
\sigma _{1}$ and $g(w)$ does not yield $\sigma _{2}$. Now 
\[
w\in \partial \left( D\right) \iff w=ti,\left\vert t\right\vert \geq \sqrt{3}%
. 
\]

Then $w_{1}=-w_{3},w_{2}=tiw_{3}$, and $p(z)=(z^{2}-w_{3}^{2})(z-itw_{3})$.
If $\func{Im}w_{3}\neq 0$, then the ratios are defined, and a simple
computation shows that 
\begin{equation}
\sigma _{1}=\dfrac{1}{3}\dfrac{it+i\sqrt{t^{2}-3}+3}{it+1},w\in \partial
\left( D\right)  \tag{(7)}
\end{equation}%
One can also compute $\sigma _{2}$ using (5), but we shall not require that
here.

\textbf{Notation: }We write $\sigma _{1}=\sigma _{1}(w)$ or $\sigma
_{2}=\sigma _{2}(w)$ if $(\sigma _{1},\sigma _{2})$ is the ratio vector of $%
p(w)=(w-w_{1})(w-w_{2})(w-w_{3})$ with $w_{1}+w_{3}=0,\func{Re}w_{1}<0<\func{%
Re}w_{3}$, and $w=\dfrac{w_{2}}{w_{3}}$

We should note here that not every $w\in D$ satisfies $w=\dfrac{w_{2}}{w_{3}}
$ for some admissible pair $(w_{2},w_{3})$. For example, $w=2$ cannot occur
since $w_{2}=2w_{3}\Rightarrow \func{Re}w_{2}>\func{Re}w_{3}$. Of course the
bounds we derive for $w\in D\cup \partial \left( D\right) $ then apply to
the subset of values of $w$ which can arise from admissible pairs. In
addition, there are admissible pairs $(w_{2},w_{3})$ such that $w\partial
\left( D\right) $, such as $w_{2}=2i,w_{3}=1$. This is not a problem since
the bounds we derive below are for $w\in D\cup \partial \left( D\right) $.
Finally, the ratios themselves are not defined when $w=1$ or $w=-1$(else the 
$w_{k}$ are not distinct). The real and imaginary parts of $f$ and of $g$
are each harmonic functions, and we want to apply the Maximum--Minimum
Principle for harmonic functions to find bounds on the real and imaginary
parts of $\sigma _{1}$ and $\sigma _{2}$. Since $D$ is unbounded, we shall
require the following special case of the Maximum--Minimum Principle for
possibly unbounded domains(see [2], page 8, Corollary 1.10]).

\textbf{Proposition 1:} Let $u$ be a real--valued harmonic function in a
domain $D$ in $R^{2}$ and suppose that 
\[
\limsup\limits_{k\rightarrow \infty }u(a_{k})\leq M 
\]%
for every sequence $\left\{ a_{k}\right\} $ in $D$ converging to a point in $%
\partial \left( D\right) $ or to $\infty $. Then $u\leq M$ on $D$.

\textbf{Remark: }As noted in [2], Proposition 1 remains valid if "$\limsup $%
" is replaced by "$\liminf $" and the inequalities are reversed.

We also need the following Local Maximum--Minimum Principle for harmonic
functions for possibly unbounded domains(see [2], page 23) to prove the
sharpness of our bounds on the real and imaginary parts of $\sigma _{1}$ and 
$\sigma _{2}$. One can prove these bounds directly, but that involves a two
variable optimization problem. Using the Maximum--Minimum Principle reduces
it to a one variable optimization problem.

\textbf{Proposition 2:} Let $u$ be a real--valued harmonic function in a
domain $D$ in $R^{2}$ and suppose that $u$ has a local maximum(or minimum)
in $D$. Then $u$ is constant.

First we require the following lemmas.

\textbf{Lemma 1}: (A) The equation $4t\sqrt{t^{2}-3}-5t^{2}+3=0$ has no real
solutions.

(B) The equation $4t\sqrt{t^{2}-3}+5t^{2}-3=0$ has no real solutions.

\textbf{Proof: }$4t\sqrt{t^{2}-3}=5t^{2}-3\Rightarrow 16t^{2}\left(
t^{2}-3\right) -(5t^{2}-3)^{2}=0\Rightarrow -9\left( t^{2}+1\right) ^{2}=0$,
which has no real solutions. That proves (A), and (B) follows in a\ similar
fashion.

\textbf{Lemma 2}: (A)The only real solution of the equation $%
t^{3}-7t-2\left( t^{2}-1\right) \sqrt{t^{2}-3}=0$ is $t=-2$.

(B) The only real solution of the equation $t^{3}-7t+2\left( t^{2}-1\right) 
\sqrt{t^{2}-3}=0$ is $t=2$.

\textbf{Proof: }$t^{3}-7t=2\left( t^{2}-1\right) \sqrt{t^{2}-3}\Rightarrow
(t^{3}-7t)^{2}-4\left( t^{2}-1\right) ^{2}(t^{2}-3)=0\Rightarrow $

$-3\left( t-2\right) \left( t+2\right) \left( t^{2}+1\right) ^{2}=0$. $t=-2$
is a solution of the given equation, but not $t=2$. That proves (A), and (B)
follows in a\ similar fashion.

\textbf{Theorem 1:} Let $p(w)=(w-w_{1})(w-w_{2})(w-w_{3}),$with $\func{Re}%
w_{1}<\func{Re}w_{2}<\func{Re}w_{3}$. Let $z_{1}$ and $z_{2}$ be the
critical points of $p$, where $z_{1}=z_{2}$ or $\func{Re}z_{1}<\func{Re}%
z_{2} $ if $z_{1}\neq z_{2}$. Let $\sigma _{1}=\dfrac{z_{1}-w_{1}}{%
w_{2}-w_{1}}$. Then

(A) $0<\func{Re}\sigma _{1}<\dfrac{2}{3}$ and the inequality is sharp in
that there are $w_{1},w_{2},$ and $w_{3}$ satisfying the hypotheses above
and such that $\func{Re}\sigma _{1}$ can be made arbitrarily close to $0$ or
arbitrarily close to $\dfrac{2}{3}$.

(B)\textit{\ }$-\dfrac{1}{3}\leq \func{Im}\sigma _{1}\leq \dfrac{1}{3}$.

(C) $\func{Im}\sigma _{1}=\dfrac{1}{3}\iff $ the roots of $p$ have the form $%
\pm i\left( z_{0}+C\right) $ and $2\left( z_{0}+C\right) $, where $\func{Im}%
z_{0}<0$, $0<\func{Re}z_{0}<-\dfrac{1}{2}\func{Im}z_{0},$ and $C$ is an
arbitrary constant.

(D) $\func{Im}\sigma _{1}=-\dfrac{1}{3}\iff $ the roots of $p$ have the form 
$\pm i\left( z_{0}+C\right) $ and $2\left( z_{0}+C\right) $, where $\func{Im}%
z_{0}>0$, $0<\func{Re}z_{0}<\dfrac{1}{2}\func{Im}z_{0},$ and $C$ is an
arbitrary constant.

(E) $\left\vert \sigma _{1}\right\vert \leq \dfrac{2}{3}$

\textbf{Proof: }While it is not necessary for $f$ to extend to be continuous
on $\partial \left( D\right) $ to apply Proposition 1, we must show that $%
0<f(w)<\dfrac{2}{3}$ for $w\in D\cup \partial \left( D\right) $ since $%
\sigma _{1}$ can arise for $w\in \partial \left( D\right) $. First we
consider the behavior of $f$ at $\infty $. $\lim\limits_{w\rightarrow \infty
}f(w)=\lim\limits_{w\rightarrow 0}f(1/w)=\dfrac{1}{3}\lim\limits_{w%
\rightarrow 0}\dfrac{1+3w\pm \sqrt{3w^{2}+1}}{w+1}=0$ or $\dfrac{2}{3}$
depending upon whether $w\rightarrow 0$ through $\func{Re}w>0$ or $\func{Re}%
w<0$. Thus by Proposition 1, $\limsup\limits_{k\rightarrow \infty }\func{Re}%
f(a_{k})\leq \dfrac{2}{3},\liminf\limits_{k\rightarrow \infty }\func{Re}%
f(a_{k})\geq 0,\limsup\limits_{k\rightarrow \infty }\func{Im}f(a_{k})\leq 
\dfrac{1}{3},$ and $\liminf\limits_{k\rightarrow \infty }\func{Im}%
f(a_{k})\geq -\dfrac{1}{3}$ for any sequence $\left\{ a_{k}\right\} $ in $D$
converging to $\infty $. We now show that $0\leq \func{Re}f\leq \dfrac{2}{3}$
and $-\dfrac{1}{3}\leq \func{Im}f\leq \dfrac{1}{3}$ as $w$ \textit{%
approaches }any point $z\in \partial \left( D\right) $. As $w$ approaches $%
z\in \partial \left( D\right) $, $\sqrt{3+w^{2}}$ approaches $\pm \sqrt{%
3-t^{2}}=\pm i\sqrt{t^{2}-3}$. Thus $\dfrac{1}{3}\dfrac{w+3-\sqrt{3+w^{2}}}{%
w+1}$ approaches $\dfrac{1}{3}\dfrac{ti+3\pm i\sqrt{t^{2}-3}}{ti+1}$. Note
that $\dfrac{1}{3}\dfrac{ti+3+i\sqrt{t^{2}-3}}{ti+1}=\sigma _{1}(w),w\in
\partial \left( D\right) $ by (7). Thus by finding the maximum and minimum
of $\func{Re}\dfrac{1}{3}\dfrac{ti+3\pm i\sqrt{t^{2}-3}}{ti+1}$ and $\func{Im%
}\dfrac{1}{3}\dfrac{ti+3\pm i\sqrt{t^{2}-3}}{ti+1},\left\vert t\right\vert
\geq \sqrt{3}$, we are finding the maximum and minimum of $\func{Re}f(w)$
and of $\func{Im}f(w)$ as $w$ approaches $\partial \left( D\right) $, and
the maximum and minimum of $\func{Re}\sigma _{1}$ and of $\func{Im}\sigma
_{1}$ for $w\in \partial \left( D\right) $. Now 
\[
\dfrac{1}{3}\dfrac{ti+3+i\sqrt{t^{2}-3}}{ti+1}=u_{1}(t)+iv_{1}(t),\dfrac{1}{3%
}\dfrac{ti+3-i\sqrt{t^{2}-3}}{ti+1}=u_{2}(t)+iv_{2}(t) 
\]%
where

\begin{equation}
u_{1}(t)=\dfrac{1}{3}\dfrac{t^{2}+3+t\sqrt{t^{2}-3}}{t^{2}+1},u_{2}(t)=%
\dfrac{1}{3}\dfrac{t^{2}+3-t\sqrt{t^{2}-3}}{t^{2}+1}  \tag{(8)}
\end{equation}%
and

\begin{equation}
v_{1}(t)=\dfrac{1}{3}\dfrac{-2t+\sqrt{t^{2}-3}}{t^{2}+1},v_{2}(t)=\dfrac{1}{3%
}\dfrac{-2t-\sqrt{t^{2}-3}}{t^{2}+1}  \tag{(9)}
\end{equation}%
$u_{1}^{\prime }(t)=-\dfrac{1}{3}\dfrac{4t\sqrt{t^{2}-3}-5t^{2}+3}{\sqrt{%
t^{2}-3}\left( t^{2}+1\right) ^{2}}$ and $u_{2}^{\prime }(t)=-\dfrac{1}{3}%
\dfrac{4t\sqrt{t^{2}-3}+5t^{2}-3}{\sqrt{t^{2}-3}\left( t^{2}+1\right) ^{2}}$%
. By Lemma 1, $u_{1}^{\prime }$ and $u_{2}^{\prime }$ have no real roots,
and hence $u_{1}$ and $u_{2}$ have no real critical points. Now $u_{1}(\sqrt{%
3})=u_{1}(-\sqrt{3})=u_{2}(\sqrt{3})=u_{2}(-\sqrt{3})=\dfrac{1}{2}$, $%
\lim\limits_{t\rightarrow \infty }u_{1}(t)=\lim\limits_{t\rightarrow -\infty
}u_{2}(t)=\dfrac{2}{3}$, and $\lim\limits_{t\rightarrow -\infty
}u_{1}(t)=\lim\limits_{t\rightarrow \infty }u_{2}(t)=0$. Thus $0\leq
u_{1}(t),u_{2}(t)\leq \dfrac{2}{3}$ for $\left\vert t\right\vert \geq \sqrt{3%
}$, which implies that $0\leq \func{Re}\dfrac{1}{3}\dfrac{ti+3\pm i\sqrt{%
t^{2}-3}}{ti+1}\leq \dfrac{2}{3}$for $\left\vert t\right\vert \geq \sqrt{3}$%
. It follows that $\limsup\limits_{k\rightarrow \infty }\func{Re}%
f(a_{k})\leq \dfrac{2}{3}$ and $\liminf\limits_{k\rightarrow \infty }\func{Re%
}f(a_{k})\geq 0$ for any sequence $\left\{ a_{k}\right\} $ in $D$ converging
to $\partial \left( D\right) $. By Proposition 1, $0\leq f(w)\leq \dfrac{2}{3%
}$ for $w\in D$. As noted above, the same proof shows that $0\leq \func{Re}%
f\leq \dfrac{2}{3}$ for $w\in \partial \left( D\right) $. By Proposition 2, $%
0<\func{Re}f<\dfrac{2}{3}$ for $w\in D$. It also follows easily that $u_{2}$
is increasing for $t\leq -\sqrt{3}$ and decreasing for $t\geq \sqrt{3}$,
which implies that $u_{2}(t)\neq 0$ and $u_{2}(t)\neq \dfrac{2}{3}$ for $%
\left\vert t\right\vert \geq \sqrt{3}$. Since $u_{2}(t)=\func{Re}%
f(w),w=ti,\left\vert t\right\vert \geq \sqrt{3},0<\func{Re}f<\dfrac{2}{3}$
for $w\in \partial \left( D\right) $. That shows that $0<\func{Re}\sigma
_{1}<\dfrac{2}{3}$. To finish the proof of part (A), if $t>\sqrt{3}$, let $%
w_{1}=-2t-i$, $w_{2}=-t+2t^{2}i$, and $w_{3}=2t+i$, while if $t<-\sqrt{3}$,
let $w_{1}=2t+i$, $w_{2}=t-2t^{2}i$, and $w_{3}=-2t-i$. In either case, $%
w=ti $ and $\func{Im}\left( 3w_{3}^{2}+w_{2}^{2}\right) =\allowbreak
12t-4t^{3}\neq 0\Rightarrow \func{Re}\sqrt{3w_{3}^{2}+w_{2}^{2}}\neq 0$.
Thus $z_{1}$ and $z_{2}$ have unequal real parts. Since $\func{Re}w_{1}<%
\func{Re}w_{2}<\func{Re}w_{3}$ as well, the ratios are defined. Above we
showed that 
\begin{equation}
\sigma _{1}(w)=u_{1}(t)+iv_{1}(t),w=it,\left\vert t\right\vert \geq \sqrt{3}
\tag{(10)}
\end{equation}%
Thus $\func{Re}\sigma _{1}(w)=$ $u_{1}(t)$. Since $\lim\limits_{t\rightarrow
\infty }u_{1}(t)=\dfrac{2}{3}$ and $\lim\limits_{t\rightarrow -\infty
}u_{2}(t)=0$, we can make $\func{Re}\sigma _{1}$ as close to $\allowbreak 0$
or $\dfrac{2}{3}$ by taking $\left\vert t\right\vert $ sufficiently large.
That finishes the proof of part (A).

To prove part (B), $v_{1}^{\prime }(t)=\allowbreak \dfrac{1}{3}\dfrac{2\sqrt{%
t^{2}-3}t^{2}-2\sqrt{t^{2}-3}-t^{3}+7t}{\sqrt{t^{2}-3}\left( t^{2}+1\right)
^{2}}$ and $v_{2}^{\prime }(t)=\dfrac{1}{3}\dfrac{2\sqrt{t^{2}-3}t^{2}-2%
\sqrt{t^{2}-3}+t^{3}-7t}{\sqrt{t^{2}-3}\left( t^{2}+1\right) ^{2}}$. By
Lemma 2, $v_{1}$ has one real critical point, $t=-2$ and $v_{2}$ has one
real critical point, $t=2$. Also, $v_{1}(\sqrt{3})=\allowbreak -\dfrac{1}{6}%
\sqrt{3}$, $v_{1}(-\sqrt{3})=\allowbreak \dfrac{1}{6}\sqrt{3}$, $v_{1}(-2)=%
\dfrac{1}{3}$, and $\lim\limits_{t\rightarrow -\infty
}v_{1}(t)=\lim\limits_{t\rightarrow \infty }v_{1}(t)=\allowbreak 0$, while $%
v_{2}(\sqrt{3})=-\dfrac{1}{6}\sqrt{3}$, $v_{2}(-\sqrt{3})=\allowbreak \dfrac{%
1}{6}\sqrt{3}$, $v_{2}(2)=\allowbreak -\dfrac{1}{3}$, and $%
\lim\limits_{t\rightarrow -\infty }v_{2}(t)=\lim\limits_{t\rightarrow \infty
}v_{2}(t)=\allowbreak 0$. Hence $-\dfrac{1}{3}\leq v_{1}(t),v_{2}(t)\leq 
\dfrac{1}{3}$ for $\left\vert t\right\vert \geq \sqrt{3}$. Arguing as
earlier, by Proposition 1 that proves part (B).

To prove (C), suppose that $\func{Im}\sigma _{1}=\dfrac{1}{3}$. If $\sigma
_{1}=\sigma _{1}(w),w\in D$, then $\func{Im}f(w)=\dfrac{1}{3}$, which cannot
happen by Proposition 2. If $\sigma _{1}=\sigma _{1}(w),w\in \partial \left(
D\right) $, then $v_{1}(t)=\allowbreak \dfrac{1}{3}$ by (7). Now it follows
easily that the only real solution of $v_{1}(t)=\allowbreak \dfrac{1}{3}$ is 
$t=-2$, and $t=-2\Rightarrow w=-2i\Rightarrow w_{3}=\dfrac{1}{2}%
iw_{2},w_{1}=-\dfrac{1}{2}iw_{2}$. The critical points of the coresponding $%
p $ are $z=\dfrac{1}{2}w_{2}$ and $z=\dfrac{1}{6}w_{2}$, which have unequal
real parts if $\func{Re}w_{2}\neq 0$. $\func{Re}w_{3}>0\Rightarrow -\dfrac{1%
}{2}\func{Im}w_{2}>0\Rightarrow \func{Im}w_{2}<0$. Also, $\func{Re}w_{2}<%
\func{Re}w_{3}\Rightarrow \func{Re}w_{2}<-\dfrac{1}{2}\func{Im}w_{2}$. If $%
\func{Re}w_{2}<0$, then $z_{1}=\dfrac{1}{2}w_{2}$ and $z_{2}=\dfrac{1}{6}%
w_{2}\Rightarrow \sigma _{1}=\dfrac{\dfrac{1}{2}w_{2}+\dfrac{1}{2}iw_{2}}{%
w_{2}+\dfrac{1}{2}iw_{2}}=\dfrac{3}{5}+\dfrac{1}{5}i\Rightarrow \func{Im}%
\sigma _{1}\neq \dfrac{1}{3}$. Letting $z_{0}=\dfrac{1}{2}w_{2}$, that
yields roots of the form $\pm iz_{0}$ and $2z_{0}$, where $\func{Re}z_{0}>0$
and $\func{Re}z_{0}<-\dfrac{1}{2}\func{Im}z_{0}$. Since any translation of $%
p $ yields the same ratios, the roots of $p$ must have the form given in
part (C). If the roots of $p$ have the form given in part (C), then $z_{1}=%
\dfrac{1}{6}w_{2}$ and $z_{2}=\dfrac{1}{2}w_{2}$, which implies that $\sigma
_{1}=\dfrac{\dfrac{1}{6}w_{2}+\dfrac{1}{2}iw_{2}}{w_{2}+\dfrac{1}{2}iw_{2}}=%
\dfrac{1}{3}+\dfrac{1}{3}i\Rightarrow \func{Im}\sigma _{1}=\dfrac{1}{3}$.
The proof of part (D) follows in a simialr fashion and we omit it.

Finally, to prove (E), note first that $f(w)=0\Rightarrow w+3-\sqrt{3+w^{2}}%
=0\Rightarrow \left( w+3\right) ^{2}-\left( 3+w^{2}\right) =\allowbreak
6w+6=0\Rightarrow w=-1$, but $f(-1)=\dfrac{1}{2}\neq 0$. Thus $f$ has no
zero in $D$ and by ([6], Theorem 13.12, page 294), $\log \left\vert
f\right\vert $ is harmonic in $D$. We shall apply Proposition 1 to $\log
\left\vert f\right\vert $. Since we showed earlier that $\lim\limits_{w%
\rightarrow \infty }f(w)=0$ or $\dfrac{2}{3}$, $\limsup\limits_{k\rightarrow
\infty }\log \left\vert f\right\vert (a_{k})\leq \log \dfrac{2}{3}$ for any
sequence $\left\{ a_{k}\right\} $ in $D$ converging to $\infty $. As $w$
approaches $\partial \left( D\right) $, $9\left\vert f(w)\right\vert ^{2}$
approaches $9\left[ \left( u_{1}(t)\right) ^{2}+\left( v_{1}(t)\right) ^{2}%
\right] $ or $9\left[ \left( u_{2}(t)\right) ^{2}+\left( v_{2}(t)\right) ^{2}%
\right] $, where $w=it,\left\vert t\right\vert \geq \sqrt{3}$. Now $%
\allowbreak 9\left[ \left( u_{1}(t)\right) ^{2}+\left( v_{1}(t)\right) ^{2}%
\right] =a(t)=\allowbreak 2\dfrac{t^{2}+3+t\sqrt{t^{2}-3}}{t^{2}+1}$, and it
follows easily that $a^{\prime }(t)=\allowbreak 2\dfrac{-4t\sqrt{t^{2}-3}%
+5t^{2}-3}{\sqrt{t^{2}-3}\left( t^{2}+1\right) ^{2}}>0$ for all $%
t,\left\vert t\right\vert \geq \sqrt{3}$. Since $\lim\limits_{t\rightarrow
\infty }a(t)=\allowbreak 4$, $a(t)<4$ for all $t,\left\vert t\right\vert
\geq \sqrt{3}$. $9\left[ \left( u_{2}(t)\right) ^{2}+\left( v_{2}(t)\right)
^{2}\right] =b(t)=\allowbreak 2\dfrac{t^{2}+3-t\sqrt{t^{2}-3}}{t^{2}+1}$. It
also follows easily that $b^{\prime }(t)=\allowbreak -2\dfrac{4t\sqrt{t^{2}-3%
}+5t^{2}-3}{\sqrt{t^{2}-3}\left( t^{2}+1\right) ^{2}}<0$ for all $%
t,\left\vert t\right\vert \geq \sqrt{3}$. Since $\lim\limits_{t\rightarrow
-\infty }b(t)=\allowbreak 4$, $b(t)<4$ for all $t,\left\vert t\right\vert
\geq \sqrt{3}$. Hence $\left\vert f(a_{k})\right\vert ^{2}\leq \dfrac{4}{9}$
for any sequence $\left\{ a_{k}\right\} $ in $D$ converging to $\partial
\left( D\right) $, which implies that $\limsup\limits_{k\rightarrow \infty
}\log \left\vert f(a_{k})\right\vert \leq \dfrac{1}{2}\log \dfrac{4}{9}$.
That proves that $\left\vert f(w)\right\vert \leq \dfrac{2}{3},w\in D$, by
Proposition 1. Note also that $9\left\vert \sigma _{1}\right\vert ^{2}=9%
\left[ \left( u_{1}(t)\right) ^{2}+\left( v_{1}(t)\right) ^{2}\right] $ for $%
w\in \partial \left( D\right) $. By what we just proved, $\left\vert \sigma
_{1}(w)\right\vert \leq \dfrac{2}{3},w\in \partial \left( D\right) $. That
finishes the proof of part (E).

\textbf{Theorem 2: }Let $p(w)=(w-w_{1})(w-w_{2})(w-w_{3}),$with $\func{Re}%
w_{1}<\func{Re}w_{2}<\func{Re}w_{3}$. Let $z_{1}$ and $z_{2}$ be the
critical points of $p$, where $z_{1}=z_{2}$ or $\func{Re}z_{1}<\func{Re}%
z_{2} $ if $z_{1}\neq z_{2}$. Let $\sigma _{2}=\dfrac{z_{2}-w_{2}}{%
w_{3}-w_{2}}$. Then

(A) $\dfrac{1}{3}<\func{Re}\sigma _{2}<1$ and the inequality is sharp in
that there are $w_{1},w_{2},$ and $w_{3}$ satisfying the hypotheses above
and such that $\func{Re}\sigma _{2}$ can be made arbitrarily close to $%
\dfrac{1}{3}$ or arbitrarily close to $1$.

(B) $-\dfrac{1}{3}\leq \func{Im}\sigma _{2}\leq \dfrac{1}{3}$

(C) $\func{Im}\sigma _{2}=\dfrac{1}{3}\iff $ the roots of $p$ have the form $%
\pm iz$ and $2z$, where $\func{Im}z<0$ and $0<\func{Re}z<-\dfrac{1}{2}\func{%
Im}z$

(D) $\func{Im}\sigma _{2}=-\dfrac{1}{3}\iff $ the roots of $p$ have the form 
$\pm iz$ and $2z$, where $\func{Im}z>0$ and $\func{Re}z>0$ and $0<\func{Re}z<%
\dfrac{1}{2}\func{Im}z$

(E) $\left\vert \sigma _{2}\right\vert \leq 1$

\textbf{Proof: }We proceed exactly as in the proof of Theorem1, working with 
$g(w)$ instead of with $f(w)$. Since\textbf{\ }$\lim\limits_{w\rightarrow
\infty }f(w)=0$ or $\dfrac{2}{3}$, by (6), $\lim\limits_{w\rightarrow \infty
}g(w)=\dfrac{1}{3}$ or $1$. As $w$ approaches $z\in \partial \left( D\right) 
$, $g(w)$ approaches $\dfrac{1}{3}\dfrac{-2ti\pm \sqrt{3-t^{2}}}{1-ti}=%
\dfrac{1}{3}\dfrac{-2ti\pm i\sqrt{t^{2}-3}}{1-ti}=1-u_{1}(t)+iv_{1}(t)$ or $%
1-u_{2}(t)+iv_{2}(t)$. Since we showed that $0<u_{1}(t)<\dfrac{2}{3},$ $%
0<u_{2}(t)<\dfrac{2}{3},-\dfrac{1}{3}\leq v_{1}(t)\leq \dfrac{1}{3},$ and $-%
\dfrac{1}{3}\leq v_{2}(t)\leq \dfrac{1}{3}$ for $\left\vert t\right\vert
\geq \sqrt{3}$, it follows immediately that $\dfrac{1}{3}<\func{Re}\sigma
_{2}<1$ and $-\dfrac{1}{3}\leq \func{Im}\sigma _{2}\leq \dfrac{1}{3}$. The
rest of parts (A) and (B) follow as in the proof of Theorem 1, parts (A) and
(B). Parts (C) and (D) also follow as in the proof of Theorem 1 parts (C)
and (D), and part (E) follows directly from Theorem 1, part (E) and (5).

\textbf{Theorem 3: }Let $p(w)=(w-w_{1})(w-w_{2})(w-w_{3}),$with $\func{Re}%
w_{1}<\func{Re}w_{2}<\func{Re}w_{3}$. Let $z_{1}$ and $z_{2}$ be the
critical points of $p$, where $z_{1}=z_{2}$ or $\func{Re}z_{1}<\func{Re}%
z_{2} $ if $z_{1}\neq z_{2}$. Let $\sigma _{1}=\dfrac{z_{1}-w_{1}}{%
w_{2}-w_{1}}$ and $\sigma _{2}=\dfrac{z_{2}-w_{2}}{w_{3}-w_{2}}$. Then $%
\func{Re}\sigma _{2}\geq \func{Re}\sigma _{1}$.

\textbf{Proof: }First, $g(w)-f(w)=\allowbreak \dfrac{1}{3}\dfrac{w^{2}+3-2%
\sqrt{3+w^{2}}}{w^{2}-1}$ is analytic in $D$ which implies that $\func{Re}%
\left( g(w)-f(w)\right) $ is a harmonic function in $D$, so we my apply
Proposition 1. Now $\lim\limits_{w\rightarrow \infty }\left(
g(w)-f(w)\right) =\allowbreak \dfrac{1}{3}\lim\limits_{w\rightarrow \infty }%
\dfrac{w^{2}+3-2\sqrt{3+w^{2}}}{w^{2}-1}=\allowbreak \dfrac{1}{3}\geq 0$.
Thus $\liminf\limits_{k\rightarrow \infty }\func{Re}\left(
f(a_{k})-g(a_{k})\right) \geq 0$ for any sequence $\left\{ a_{k}\right\} $
in $D$ converging to $\infty $. Also, as $w\rightarrow \partial \left(
D\right) $, $g(w)-f(w)\rightarrow \dfrac{1}{3}\dfrac{-t^{2}+3\pm 2i\sqrt{%
t^{2}-3}}{-t^{2}-1}=\dfrac{1}{3}\dfrac{t^{2}-3\pm 2i\sqrt{t^{2}-3}}{t^{2}+1}$%
. Then $\func{Re}\left( g(w)-f(w)\right) \rightarrow \dfrac{1}{3}\dfrac{%
t^{2}-3}{t^{2}+1}\geq 0$ since $t^{2}\geq 3$. It follows that $%
\liminf\limits_{k\rightarrow \infty }\func{Re}\left(
f(a_{k})-g(a_{k})\right) \geq 0$ for any sequence $\left\{ a_{k}\right\} $
in $D$ converging to $\partial \left( D\right) $. By Proposition 1, $\func{Re%
}\left( \sigma _{2}(w)-\sigma _{1}(w)\right) \geq 0$, $w\in D$. For $w\in
\partial \left( D\right) $, by (7) and (5), $\sigma _{2}=\allowbreak -i%
\dfrac{it+1}{2t-\sqrt{t^{2}-3}}$, which implies that 
\begin{equation}
\sigma _{2}-\sigma _{1}=\allowbreak \dfrac{1}{3}\dfrac{-t^{2}+3+2i\sqrt{%
t^{2}-3}}{-t^{2}-1}.  \tag{(11)}
\end{equation}%
By what was just proved, $\func{Re}\left( \sigma _{2}(w)-\sigma
_{1}(w)\right) \geq 0$, $w\in \partial \left( D\right) $.

\textbf{Note:} The example below shows that is possible to have $\func{Re}%
\sigma _{2}=\func{Re}\sigma _{1}$. In fact, below we have $\sigma
_{1}=\sigma _{2}$.

\textbf{Example:} Let $w_{1}=-1$, $w_{2}=\sqrt{3}i$, $w_{3}=1$, which
implies that $w=\sqrt{3}i$ and the $\left\{ w_{k}\right\} $ are the vertices
of an equilateral triangle. Then $z_{1}=z_{2}=\dfrac{1}{\sqrt{3}}i$ and $%
\sigma _{1}=\sigma _{2}=\dfrac{1}{2}-\dfrac{1}{6}i\sqrt{3}$.

It is natural to ask whether the example above gives essentially the only
case when $\sigma _{1}=\sigma _{2}$.

\textbf{Theorem 4: }$\sigma _{1}=\sigma _{2}\iff w_{1},w_{2},w_{3}$ are the
vertices of an equilateral triangle which contains no vertical line segment.

\textbf{Proof: }$w=\pm 1\Rightarrow w_{2}=w_{3}$ or $w_{2}=w_{1}$, in which
case the ratios are not defined. Thus we may assume that $w\neq \pm 1$. For $%
w\in D$, $\sigma _{1}(w)=\sigma _{2}(w)\iff f(w)-g(w)=-\dfrac{1}{3}\dfrac{%
w^{2}-2\sqrt{3+w^{2}}+3}{\left( -1+w\right) \left( w+1\right) }=0\iff $

$w^{2}-2\sqrt{3+w^{2}}+3=0\iff w=\pm i\sqrt{3}\iff \left\{
w_{1},w_{2},w_{3}\right\} =\left\{ -w_{3},\pm \sqrt{3}iw_{3},w_{3}\right\} $%
, which are easily seen to be the vertices of an equilateral triangle. For $%
w\in \partial \left( D\right) $, by (11), $\sigma _{1}(w)=\sigma _{2}(w)\iff 
$

$\allowbreak \dfrac{1}{3}\dfrac{-t^{2}+3+2i\sqrt{t^{2}-3}}{-t^{2}-1}%
=0,w=ti,\left\vert t\right\vert \geq \sqrt{3}$. That yields $t=\pm \sqrt{3}$%
, which gives $w=\pm i\sqrt{3}$ as above. We can also assume that the
triangle formed by $w_{1},w_{2},w_{3}$ contains no vertcial line segment,
since the ratios are not defined in that case either.

\textbf{Theorem 5: }$\sigma _{1}$ or $\sigma _{2}$ are real if and only if $%
w_{1},w_{2},$ and $w_{3}$ are collinear.

\textbf{Proof: }Suppose first that $\sigma _{1}(w)$ is real, $w\in D$. Then $%
f(w)$ is real, or $w+3-\sqrt{3+w^{2}}=k(w+1),k\in \Re $, which implies,
after some simplification, that $(k^{2}-2k)w^{2}+2(1-k)(3-k)w+k^{2}-6k+6=0$.
The discriminant of this quadratic equation is $%
4(1-k)^{2}(3-k)^{2}-4(k^{2}-2k)(k^{2}-6k+6)=4\left( 2k-3\right) ^{2}\geq 0$
since $k\in \Re $. Hence $w$ is real. Now if $w$ is real, then for the
ratios to exist, $w=\dfrac{w_{2}}{w_{3}}$ must be a positive real number,
which we again denote by $k$. But then $w_{1}=-kw_{2}$ and $w_{3}=kw_{2}$.
It is then easy to show that the set of points $\left\{
-kw_{2},w_{2},kw_{2}\right\} $ must be collinear. If $\sigma _{1}(w)$ is
real, $w\in \partial \left( D\right) $, then by (10) and (9), $-2t+\sqrt{%
t^{2}-3}=0$, which has no real solutions. If $\sigma _{2}$ is real, we can
proceed in the same fashion, or just use (5) to show that $\sigma _{1}$ is
real.

\textbf{Remark: }One can easily extend the definition of complex ratios
given in this paper to functions of the form $%
p(z)=(z-w_{1})^{m_{1}}(z-w_{2})^{m_{2}}(z-w_{3})^{m_{3}}$, where $%
m_{1},m_{2},$ and $m_{3}$ are given positive real numbers. This is discussed
in [4] for all real $w_{j}$.

\section{References}

(1) Peter Andrews, Where not to find the critical points of a
polynomial-variation on a Putnam theme, Amer. Math. Monthly 102(1995)
155--158.

(2) Sheldon Axler, Paul Bourdon, and Wade Ramey, Harmonic Function Theory,
2001, Springer Verlag, New York, Inc.

(3) Alan Horwitz, On the Ratio Vectors of Polynomials, Journal of
Mathematical Analysis and Applications 205(1997), 568-576.

(4) Alan Horwitz, Ratio vectors of polynomial--like functions, preprint.

(5) Gideon Peyser, On the roots of the derivative of a polynomial with real
roots, Amer. Math. Monthly 74(1967), 1102--1104.

(6) Walter Rudin, Real and Complex Analysis, 2nd ed., 1974, McGraw--Hill.

\end{document}